# On Self-Predicative Universals in Category Theory

David Ellerman
Philosophy Department,
University of California/Riverside

Abstract
1. This paper shows how the universals of category theory in mathematics provide a model (in the Platonic Heaven of mathematics) for the self-predicative strand of Plato's Theory of Forms as well as for the idea of a "concrete universal" in Hegel and similar ideas of paradigmatic exemplars in ordinary thought.
2. The paper also shows how the always-self-predicative universals of category theory provide the "opposite bookend" to the never-self-predicative universals of iterative set theory and thus that the paradoxes arose from having one theory (e.g., Frege's Paradise) where universals could be either self-predicative or non-self-predicative (instead of being always one or the other).
3. Moreover the paper considers one of the most important examples of self-predicative universals in pure mathematics, namely adjoint functors or adjunctions. It gives a parsing of adjunctions into two halves (*left* and *right semi-adjunctions*) using the heterodox notion of heteromorphisms, and then shows that the parts can be recombined in a new way to define the cognate-to-adjoints notion of a *brain functor* that provides an abstract conceptual model of a brain.
4. Finally the paper argues that at least one way category theory has foundational relevance is that it isolates the universal concepts and structures that are important throughout mathematics.

**Introduction: "Bad Platonic Metaphysics"**

Consider the following example of "bad metaphysics."

> Given all the entities that have a certain property, there is one entity among them that exemplifies the property in a universal, paradigmatic, archetypical, ideal, essential, or canonical way. It is called the "self-predicative universal." There is a relationship of "participation" or "resemblance" so that all and only the entities that have the property (including itself) "participate in" or "resemble" that perfect example, the self-predicative universal.

To the modern ear, all this sounds like the worst sort of "bad Platonic metaphysics." Yet there is a mathematical theory developed within the last seventy years, category theory [MacLane 1971], that provides precisely that treatment of self-predicative universals within mathematics.

A simple example using sets will illustrate the points. Given two sets a and b, consider the property of sets: $F(x) \equiv$ "x is contained in a and is contained in b." In other words, the property is the property of being both a subset of a and a subset of b. In this example, the *participation* relation is the subset inclusion relation. There is a set, namely the intersection or meet of a and b, denoted a∩b, that has the property (so it is a "concrete" instance of the property), and it is universal in the sense that any other set has the property if and only if it participates in (i.e., is included in) the universal example:

  self-predication:  $F(a \cap b)$, i.e., $a \cap b \subseteq a$ and $a \cap b \subseteq b$, and
  universality:  x participates in a∩b if and only if $F(x)$, i.e., $x \subseteq a \cap b$ if and only if
          $x \subseteq a$ and $x \subseteq b$.

This example of a self-predicative universal is quite simple, but all this "bad metaphysical talk" has highly developed and precise models in category theory.

- This interpretation of the universals of category theory as self-predicative universals is the first point of this paper (written for a non-mathematical philosophical audience).
- In terms of the old theme of universals in philosophy, we show how the self-predicative universals of category theory provide a rigorous model (in the "Platonic Heaven" of mathematics) for the self-predicative strand in Plato's thought and the "concrete universal" synthesis in Hegel's thought as well as for the common "Form" of thought that considers a paradigmatic, canonical, iconic, archetypical, or quintessential exemplar of some property.
- We will also locate the always-self-predicative universals of category theory as the opposite bookend to the never-self-predicative sets of iterative set theory and see that the paradoxes arose from having one theory where universals could be either self-predicative or non-self-predicative (instead of being always one or always the other).
- Also we will consider one of the most important uses of self-predicative universals in pure mathematics, namely adjoint functors or adjunctions. We give a parsing of adjunctions into two halves (*left* and *right semi-adjunctions*) using the heterodox notion of heteromorphisms, and then show that the parts can be recombined in a new way to define the cognate-to-adjoints notion of a *brain functor* that abstractly models a brain.
- Finally the paper argues that at least one way category theory has foundational relevance is that it isolates the universal concepts and structures that are important throughout mathematics.

## Universals: Abstract or Concrete

### A Theory of Universals

In Plato's Theory of Ideas or Forms (ειδη), a property F has an entity associated with it, the *universal* $u_F$, which uniquely represents the property. An object x has the property F, i.e., F(x), if and only if (iff) the object x *participates* in the universal $u_F$. Let μ (from μεθεξις or methexis) represent the participation relation so

"x μ $u_F$" reads as "x participates in $u_F$".

Given a relation μ, an entity $u_F$ is said to be *a universal for the property F* (with respect to μ) if it satisfies the following *universality condition*:

for any x, x μ $u_F$ if and only if F(x).

A universal representing a property should be in some sense unique. Hence there should be an equivalence relation (≈) so that universals satisfy a *uniqueness condition:*

if $u_F$ and $u_F'$ are universals for the same F, then $u_F ≈ u_F'$.

A mathematical theory is said to be a *theory of universals* if it contains a binary relation μ and an equivalence relation ≈ so that with certain properties F there are associated entities $u_F$ satisfying the following conditions:

(I) *Universality*: for any x, x μ $u_F$ iff F(x), and
(II) *Uniqueness*: if $u_F$ and $u_F'$ are universals for the same F [i.e., satisfy (I) above], then $u_F ≈ u_F'$.

A universal $u_F$ is said to be *abstract* if it does not participate in itself, i.e., ¬($u_F$ μ $u_F$). A universal $u_F$ is *self-predicative* or *concrete* if it participates in itself, i.e., $u_F$ μ $u_F$.

### Set Theory as The Theory of Abstract Universals

There is a modern mathematical theory that readily qualifies as a theory of universals, namely (naïve) set theory. The universal representing a property F is the set of all elements with the property:

$$u_F = \{ x \mid F(x) \}.$$

The participation relation is the set membership relation usually represented by ∈. The universality condition in (naïve) set theory is the (naïve) *comprehension axiom*: there is a set y such that for any x, x∈y iff F(x). Set theory also has an *extensionality axiom*, which states that two sets with the same members are identical:

for all x, (x ∈ y iff x ∈ y') implies y = y'.

Thus if y and y' both satisfy the comprehension axiom scheme for the same F then y and y' have the same members so y = y'. Hence in set theory the uniqueness condition on universals is satisfied with the equivalence relation (≈) as equality (=) between sets. Thus naive set theory qualifies as a theory of universals.

The hope that naive set theory would provide a *general* theory of universals proved to be unfounded. The naïve comprehension axiom lead to inconsistency for such properties as

F(x) ≡ "x is not a member of x" ≡ x∉x

If R is the universal for that property, i.e., R is the set of all sets which are not members of themselves, the naive comprehension axiom yields a contradiction.

R∈R iff R∉R.
Russell's Paradox

The characteristic feature of Russell's Paradox and the other set theoretical paradoxes is the *negated* self-reference wherein the universal is allowed to qualify for the negated property represented by the universal, e.g., the Russell set R is allowed to be one of the x's in the universality relation: $x \in R$ iff $x \notin x$. The set-theoretic formulation of the paradox was not essential. Russell himself expounded the paradox in term of the property of predicates that they are *not* self-predicative, i.e., the "predicates which are not predicable of themselves" [Russell, 2010 (1903), 80]

There are several ways to restrict the naive comprehension axiom to defeat the set theoretical paradoxes, e.g., as in Russell's type theory, Zermelo-Fraenkel set theory, or von Neumann-Bernays set theory. The various restrictions are based on an iterative concept of set [Boolos 1971] which forces a set y to be more "abstract", e.g., of higher type or rank, than the elements $x \in y$. As Russell himself put it:

> It will now be necessary to distinguish (1) terms, (2) classes, (3) classes of classes, and so on ad infinitum; we shall have to hold that no member of one set is a member of any other set, and that $x \in u$ requires that x should be of a set of a degree lower by one than the set to which u belongs. Thus $x \in x$ will become a meaningless proposition; and in this way the contradiction is avoided. [Russell, 2010 (1903), 527]

Thus the universals provided by the various set theories are "abstract" universals in the technical sense that they are relatively more abstract (i.e., of higher type or rank) than the objects having the property represented by the universal.[1] Sets may not be members of themselves.[2]

With the modifications to avoid the paradoxes, a set theory still qualifies as a theory of universals. The membership relation is the participation relation so that for suitably restricted predicates, there exists a set satisfying the universality condition. Set equality serves as the equivalence relation in the uniqueness conditions. But set theory cannot qualify as a *general* theory of universals. The paradox-induced modifications turn the various set theories into theories of *abstract* (always non-self-predicative) universals since they prohibit the self-membership of sets. That clears the ground for another theory of self-predicative universals.

**Self-predicative or Concrete Universals**

Philosophy has long contemplated another type of universal, variously called a *self-predicative*, *self-participating*, or *concrete universal.* Indeed, it is a common Form of thought.

---

[1] As was correctly pointed out by Jean-Pierre Marquis, all mathematical concepts such as the "abstract universals" of set theory and the "concrete universals" of category theory [Ellerman 1988, 1995] "are abstract" [Marquis, 2006, 183, fn. 2; see also Marquis 2000] in the more usual sense. This use of the adjective "abstract" in "abstract universal" is used as a technical term to mean relatively abstract in the sense of higher type or rank in some iterative concept of sets than the instances of the property represented by the universal. And the adjective "concrete" in "concrete universal" is used in the technical sense to mean that the universal is "One among the many" instances of the property instead of "One over the many" (to borrow some of the language of the Third Man Argument considered later). Although "concrete universal" is a traditional philosophical (Hegel) and literary concept [Wimsatt 1947], we will mostly use the perhaps less-confusing adjective "self-predicative" [Malcolm 1991] or even "self-participating" interchangeably in case the adjective "concrete" is thought to refer to concrete entities like tables and chairs.

[2] Quine's system ML [1955b] allows "$V \in V$" for the universal class V, but no standard model of ML has ever been found where "$\in$" is interpreted as set membership [Hatcher 1982, Chapter 7].

The intuitive idea of a self-participating universal for a property is that it is an object that has the property and has it in such a universal sense that all other objects with the property resemble or participate in that paradigmatic, archetypal, canonical, iconic, ideal, essential, or quintessential exemplar. Such a universal $u_F$ for a property F is *self-predicative* in the sense that it has the property itself, i.e., $F(u_F)$. It is *universal* in the intuitive sense that it represents F-ness is such a perfect and exemplary manner that any object resembles or participates in the universal $u_F$ if and only if it has the property F.

The intuitive notion of a concrete universal occurs in ordinary thought as in the "all-American boy" or any case of a quintessential iconic example (e.g., Sophia Loren as "the" Italian women or Michelangelo's David or da Vinci's Mona Lisa as "the" exemplars for other properties), or when resemblance to an "defining" example becomes an adjective like "Lincolnesqe." In Greek-inspired Christian theology, there is the "Word made Flesh" [Miles 2005] together with *imitatio Christi* as the participation or resemblence relation to that concrete universal. The idea of the concrete universal is often associated with Hegel [Stern 2007] where it was one way to think about the synthesis between an abstract universal thesis and the antithesis of diverse particulars. One sensible application by Hegel was in the arts and literature to explicate the old idea that great art uses a concrete instance to universally exemplify certain human conditions (e.g., Shakespeare's Romeo and Juliet as "the" romantic tragedy) [Desmond 1986].

The notion of a self-predicative universal goes back to Plato's Theory of Forms [Vlastos 1978, 1981,1995; Malcolm 1991]. Plato's forms are often considered to be abstract or non-self-predicative universals quite distinct from and "above" the instances. In the words of one Plato scholar, "not even God can scratch Doghood behind the Ears" [Allen 1960]. But Plato did give examples of self-participation or self-predication, e.g., that Justice is just [Protagoras 330]. Moreover, Plato used expressions that indicated self-predication of universals.

> But Plato also used language which suggests not only that the Forms exist separately ($\chi\omega\rho\iota\sigma\tau\alpha$) from all the particulars, but also that each Form is a peculiarly accurate or good particular of its own kind, i.e., the standard particular of the kind in question or the model ($\pi\alpha\rho\alpha\delta\epsilon\iota\gamma\mu\alpha$) to which other particulars approximate. [Kneale and Kneale 1962, 19]

But many scholars regard the notion of a Form as *paradeigma* or self-predicative universal as an error.

> For general characters are not characterized by themselves: humanity is not human. The mistake is encouraged by the fact that in Greek the same phrase may signify both the concrete and the abstract, e.g. $\tau o\ \lambda\epsilon\upsilon\kappa o\nu$ (literally "the white") both "the white thing" and "whiteness", so that it is doubtful whether $\alpha\upsilon\tau o\ \tau o\ \lambda\epsilon\upsilon\kappa o\nu$ (literally "the white itself") means "the superlatively white thing" or "whiteness in abstraction". [Kneale and Kneale 1962, 19-20]

Thus some Platonic language is ambivalent between interpreting a form as a concrete universal ("the superlatively white thing") and an abstract universal ("whiteness in abstraction").

The literature on Plato has reached no resolution on the question of self-predication. Scholarship has left Plato on both sides of the fence; many universals are not self-predicative but some are. It is fitting that Plato should exhibit this ambivalence since the self-predication issue has only come to a head in the 20[th] century with the set theoretical antinomies. Set theory had to be reconstructed as a theory of universals that were rigidly non-self-predicative.

The reconstruction of set theory as the theory of *never*-self-predicative universals cleared the ground for a separate theory of universals that are *always* self-predicative. Such a theory of self-predicative universals would realize the self-predicative strand of Plato's Theory of Forms. The antinomies show that there cannot be one theory of universals (e.g., Frege's Paradise) that could be self-predicative or non-self-predicative because one could then consider the universal for all those universals that are non-self-predicative.

A theory of (always) self-predicative universals would have an appropriate participation relation $\mu$ so that for certain properties F, there are entities $u_F$ satisfying the *universality condition*:

$$\text{for any x, x } \mu \text{ } u_F \text{ if and only if F(x).}$$

The universality condition and $F(u_F)$ imply that $u_F$ is a *concrete* universal in the previously defined sense of being self-predicative, $u_F \mu u_F$. A theory of self-predicative universals would also have to have an equivalence relation so the self-predicative universals for the same property would be *the* universal up to that equivalence relation.

Is there a precise mathematical theory of the common Form of thought, the self-predicative universal? Our claim is that category theory is precisely that theory where the self-predicative universals are the universal constructions, usually as universal morphisms or universal arrows.

> Universal constructions appear throughout mathematics in various guises – as universal arrows to a given functor, as universal arrows from a given functor, or as universal elements of a set-valued functor. [MacLane, 1971, 55]

To keep matters simple and intuitive for the non-mathematician, all our examples will use one of the simplest examples of categories, namely partially ordered sets.[3] Consider the universe of subsets or power set $\wp(U)$ of a set U with the inclusion relation $\subseteq$ as the partial ordering relation. Given sets a and b, consider the property

$$G(x) \equiv a \subseteq x \text{ \& } b \subseteq x.$$

The participation relation is set inclusion $\subseteq$ and the union $a \cup b$ is the universal $u_F$ for this property G(x). The universality relation states that the union is the least upper bound of a and b in the inclusion ordering:

$$\text{for any x, } a \cup b \subseteq x \text{ iff } a \subseteq x \text{ \& } b \subseteq x.$$

The universal has the property it represents, i.e., $a \subseteq a \cup b$ & $b \subseteq a \cup b$, so it is a self-predicative or concrete universal. Two self-predicative universals for the same property must participate in each other. In partially ordered sets, the antisymmetry condition, $y \subseteq y'$ & $y' \subseteq y$ implies $y = y'$, means that equality can serve as the equivalence relation in the uniqueness condition for universals in a partial order.

---

[3] A binary relation $\leq$ on U is a *partial order* if for all u,u',u''∈U, it is reflexive ($u \leq u$), transitive ($u \leq u'$ and $u' \leq u''$ imply $u \leq u''$), and anti-symmetric ($u \leq u'$ and $u' \leq u$ imply $u = u'$). For less trivial examples with more of a category-theoretic flavor, see Ellerman 1988.

### Self-predicative Universals in more general categories

For the self-predicative universals of category theory,[4] the *participation relation* is the *uniquely-factors-through* relation. It can always be formulated in a suitable category as:

"x μ $u_F$" means "there exists a unique arrow x→ $u_F$".

Then x is said to *uniquely factor through* $u_F$, and the arrow x→ $u_F$ is the unique factor or participation morphism. In the universality condition, for any x, x μ $u_F$ if and only if F(x), the existence of the identity arrow $u_F$ → $u_F$ is the self-participation of the self-predicative universal that corresponds with F($u_F$), the self-predication of the property to $u_F$. In category theory, the equivalence relation used in the uniqueness condition is the isomorphism ($\cong$).[5]

It is sometimes convenient to "turn the arrows around" and use the dual definition where "x μ $u_G$" means "there exists a unique arrow $u_G$ →x" that can also be viewed as the original definition stated in the dual or opposite category. The above treatment of the intersection a∩b and the union a∪b are dual to one another. If we think of the "uniquely-factoring-through" arrows as *transferring* the property from the universal to the instances, then the transferring may go along the direction of the arrow—so the property may be said to be *transmitted* to the instance—or against the direction of the arrow—so the property may be said to be *reflected* to the instance.

Category theory as the theory of self-predicative universals has quite a different flavor from set theory, the theory of abstract universals. Given an appropriately delimited collection of all the elements with a property, set theory can postulate a more abstract entity, the set of those elements, to be the universal. But category theory cannot postulate its universals because those universals are self-predicative, i.e., are the "One among the many." Category theory must find its universals, if at all, among the entities with the property.

The "mistake" in the set-theoretic paradoxes and similar self-referential antimonies is often taken to be the self-reference.

> In all the above contradictions (which are merely selections from an indefinite number) there is a common characteristic, which we may describe as self-reference or reflexiveness. [Whitehead and Russell 1997 (1910), 61]

---

[4] In the general case, a category may be defined as follows [e.g., MacLane and Birkhoff 1988 or MacLane 1971]:
   A *category* C consists of
   (a) a set of *objects* a, b, c, ...,
   (b) for each pair of objects a,b, a set $hom_C$(a,b) = C(a,b) whose elements are represented as *arrows* or *morphsims* f: a → b,
   (c) for any f ∈ $hom_C$(a,b) and g ∈ $hom_C$(b,c), there is the *composition* gf: a→b→c in $hom_C$(a,c),
   (d) composition of arrows is an associative operation, and
   (e) for each object a, there is an arrow $1_a$ ∈ $hom_C$(a,a), called the *identity* of a, such that for any f: a→b and g:c→a, f$1_a$ = f and $1_a$ g = g.
An arrow f:a→b is an *isomorphism*, a $\cong$ b, if there is an arrow g:b→a such that fg = $1_b$ and gf = $1_a$. A *functor* is a map from one category to another that preserves composition and identities.

[5] Thus it must be verified that two concrete universals for the same property are isomorphic. By the universality condition, two concrete universals u and u' for the same property must participate in each other. Let f:u'→u and g:u→u' be the unique arrows given by the mutual participation. Then by composition gf:u'→u' is the unique arrow u'→u' but $1_{u'}$ is another such arrow so by uniqueness, gf = $1_{u'}$. Similarly, fg:u→u is the unique self-participation arrow for u so fg = $1_u$. Thus mutual participation of u and u' implies the isomorphism u $\cong$ u'.

The iterative notion of a set requires the universal for a property to be of higher type or rank than the instances so that "x∈x will become a meaningless proposition; and in this way the contradiction is avoided." To avoid the paradoxes, Whitehead and Russell postulated the *vicious circle principle*: "Whatever involves *all* of a collection must not be one of the collection." [1997, 37] But a self-predicative universal of category theory is the "One among the many" that transfers the property to all the instances so it is "impredicative" or self-predicative in violation of the vicious circle principle. Indeed, the universals of category theory are always self-predicative[6] via the identity morphisms so the question arises of how category theory avoids similar paradoxes.

All morphisms can be seen as "uniquely factoring through" themselves by the identity morphism (at either the head or tail of the arrow)—so the construction of something like "a universal morphism for all those morphisms that don't factor through themselves" would always come up empty. Abstractly put, there can be no self-predicative universal for the property of not being self-predicative—since the universal needs to have the property that is "transferred" to the instances by their "participation" in the universal and that particular negative property would *always* be defeated by the universal's identity morphism.[7] Thus the problem with the paradoxes was not the self-predication *per se* but the negated self-predication, and that is defeated in category theory by the universals being always self-predicative (by the identity morphisms). The "circle" or self-reference is not the problem if all the circles are required to be "virtuous" so that a "vicious" circle cannot arise.

**The Third Man Argument**

Much of the modern Platonic literature on self-participation and self-predication [e.g., Malcolm 1991] stems from the work of Geach [1956] and Vlastos [1981] on the Third Man argument. The name derives from Aristotle (who denied self-predication), but the argument occurs in the dialogues.

> But now take largeness itself and the other things which are large. Suppose you look at all these in the same way in your mind's eye, will not yet another unity make its appearance—a largeness by virtue of which they all appear large?
> So it would seem.
> If so, a second form of largeness will present itself, over and above largeness itself and the things that share in it, and again, covering all these, yet another, which will make all of them large. So each of your forms will no longer be one, but an indefinite number. [Parmenides, 132]

If a form is self-predicative, the participation relation can be interpreted as "resemblance." An instance has the property F because it resembles the paradigmatic example of

---

[6] The connection between the self-predicative or "impredicative" definitions (which caused the problems in naïve set theory) and the self-predicative universals of category theory has not escaped the attention of category theorists. For instance, Michael Makkai notes that the "Peano system" of natural numbers is the self-predicative universal for the property of being a "pre-Peano system": "we can say that a Peano system is distinguished among pre-Peano systems by the fact that it has *exactly one morphism to* any pre-Peano system. (An 'impredicative' definition if there ever was one!)" [Makkai 1999, 52] See MacLane and Birkhoff for a full explanation of that "Peano-Lawvere Axiom" [1967, 67] characterizing the Natural Numbers as the self-predicative universal for counting systems.

[7] From the purely syntactic viewpoint, it is "as if" all sets automatically had self-membership x∈x so the paradox-causing *negation* x∉x could never apply.

F-ness. But then, the Third Man argument contends, the common property shared by Largeness and other large things gives rise to a "One over the many", a form Largeness* such that Largeness and the large things share the common property by virtue of resembling Largeness*. And the argument repeats itself giving rise to an infinite regress of forms. A key part of the Third Man argument is what Vlastos calls the *Non-Identity thesis*:

> NI If anything has a given character by participating in a Form, it is not identical with that Form. [Vlastos 1981, 351]

It implies that Largeness* is not identical with Largeness.

P. T. Geach [1956] has developed a self-predicative interpretation of Forms as standards or norms, an idea he attributes to Wittgenstein. A stick is a meter long because it resembles, lengthwise, the standard meter measure. Geach avoids the Third Man regress with the exceptionalist device of holding the Form "separate" from the many so they could not be grouped together to give rise to a new "One over the many." Geach aptly notes the analogy with Frege's ad hoc and unsuccessful attempt to avoid the Russell-type paradoxes by allowing a set of all and only the sets which are not members of themselves—except for that set itself [Quine 1955a; Geach 1980].

Category theory provides a mathematical model for the Third Man argument, and it shows how to avoid the regress. The category-theoretic model shows that the flaw in the Third Man argument lies not in self-predication but in the Non-Identity thesis [Vlastos 1954, 326-329]. "The One" is not necessarily "over the many"; it can be (isomorphic to) one among the many. In mathematical terms, a colimit or limit can "take on" one of the elements in the diagram. In the special case of sets ordered by inclusion, the union or intersection of a collection of sets is not necessarily distinct from the sets in the collection; the "One" could be one *among* the many.

For example, let $A = \cup \{A_\beta\}$ be the One formed as the union of a collection of many sets $\{A_\beta\}$. Then add A to the collection and form the new One* as
$$A^* = \cup\{A_\beta\} \cup A.$$
This operation leads to no Third Man regress since $A^* = A$.

Whitehead described European philosophy as a series of footnotes to Plato, and the Theory of Forms was central to Plato's thought. We have seen that the self-predicative universals of category theory provide a rigorous mathematical model for the self-predicative strain in Plato's Theory of Forms and for the intuitive notion of a concrete universal elsewhere in philosophy, literature, and ordinary thought.

**Adjoint Functors**

### The Ubiquity and Importance of Adjoints

One of the most important and beautiful notions in category theory is the notion of a pair of adjoint functors. The developers of category theory, Saunders MacLane and Samuel Eilenberg, famously said that categories were defined in order to define functors, and functors were defined in order to define natural transformations [Eilenberg and MacLane 1945]. Adjoints were defined more than a decade later by Daniel Kan [1958] but the realization of their ubiquity ("Adjoint functors arise everywhere" [MacLane 1971, v]) and their foundational importance has steadily increased over time [Lawvere 1969; Lambek 1981]. Now it would perhaps not be too much of an exaggeration to see categories, functors, and natural transformations as the prelude to defining adjoint functors. The notion of adjoint functors (and the constituent semi-adjunctions defined below) includes all the instances of self-predicative universal mapping properties discussed above. As Steven Awodey put it:

> The notion of adjoint functor applies everything that we have learned up to now to unify and subsume all the different universal mapping properties that we have encountered, from free groups to limits to exponentials. But more importantly, it also captures an important mathematical phenomenon that is invisible without the lens of category theory. Indeed, I will make the admittedly provocative claim that adjointness is a concept of fundamental logical and mathematical importance that is not captured elsewhere in mathematics. [Awodey 2006, 179]

Other category theorists have given similar testimonials.

> To some, including this writer, adjunction is the most important concept in category theory. [Wood 2004, 6]

> The isolation and explication of the notion of *adjointness* is perhaps the most profound contribution that category theory has made to the history of general mathematical ideas." [Goldblatt 2006, 438]

> Nowadays, every user of category theory agrees that [adjunction] is the concept which justifies the fundamental position of the subject in mathematics. [Taylor 1999, 367]

**Adjoints and universals**

How do the ubiquitous and important adjoint functors relate to our theme of self-predicative universals? MacLane and Birkhoff succinctly state the idea of the self-predicative universals of category theory and note that adjunctions can be analyzed in terms of those universals.

> The construction of a new algebraic object will often solve a specific problem in a universal way, in the sense that every other solution of the given problem is obtained from this one by a unique homomorphism. The basic idea of an adjoint functor arises from the analysis of such universals. [MacLane and Birkhoff 1988, v]

We will use a specific novel treatment of adjunctions [Ellerman 2006] that shows they arise by gluing together in a certain way two universal constructions or self-predicative universals ("semi-adjunctions"). But for illustration, we will stay within the methodological restriction of using examples from partial orders (where adjunctions are called "Galois connections").

We have been working within the inclusion partial order on the set of subsets $\wp(U)$ of a universe set U. Consider the set of all ordered pairs of subsets <a,b> from the Cartesian product $\wp(U) \times \wp(U)$ where the partial order (using the same symbol $\subseteq$) is defined by pairwise inclusion. That is, given the two ordered pairs <a', b'> and <a,b>, we define[8]

---

[8] The "iff" or "if and only if" is replaced by a natural isomorphism in the general case of an adjunction or a semi-adjunction (defined later). Our focus here is on the universal constructions isolated and defined by category theory, but the notion of a natural transformation is also a key part in the characterization of the universal constructions.

$$\langle a',b'\rangle \subseteq \langle a,b\rangle \text{ if } a' \subseteq a \text{ and } b' \subseteq b.$$

Order-preserving maps can be defined each way between these two partial orders. From $\wp(U)$ to $\wp(U)\times\wp(U)$, there is the diagonal map $\Delta(x) = \langle x,x\rangle$, and from $\wp(U)\times\wp(U)$ to $\wp(U)$, there is the meet map $\cap(\langle a,b\rangle) = a \cap b$. Consider now the following "*adjointness relation*" between the two partial orders:

$$\Delta(c) \subseteq \langle a,b\rangle \text{ iff } c \subseteq \cap(\langle a,b\rangle)$$
Adjointness Equivalence

for sets a, b, and c in $\wp(U)$. It has a certain symmetry that can be exploited. If we fix $\langle a,b\rangle$, then we have the previous universality condition for the meet of a and b: for any c in $\wp(U)$,

$$c \subseteq a \cap b \text{ iff } \Delta(c) \subseteq \langle a,b\rangle.$$
Universality Condition for Meet of Sets a and b

The defining property on elements c of $\wp(U)$ is that $\Delta(c) \subseteq \langle a,b\rangle$ (just a fancy way of saying that "c is a subset of both the given a and b").

But using the symmetry, we could fix c and have another universality condition using the reverse inclusion in $\wp(U)\times\wp(U)$ as the participation relation: for any $\langle a,b\rangle$ in $\wp(U)\times\wp(U)$,

$$\langle a,b\rangle \supseteq \Delta(c) \text{ iff } c \subseteq a \cap b$$
Universality Condition for $\Delta(c)$.[9]

Here the defining property on elements $\langle a,b\rangle$ of $\wp(U)\times\wp(U)$ is that "the meet of a and b is a superset of the given set c." The self-predicative universal for that property is the image of c under the diagonal map $\Delta(c) = \langle c,c\rangle$, just as the self-predicative universal for the other property defined given $\langle a,b\rangle$ was the image of $\langle a,b\rangle$ under the meet map $\cap(\langle a,b\rangle) = a \cap b$.

Thus in this adjoint situation between the two categories $\wp(U)$ and $\wp(U)\times\wp(U)$, we have a pair of maps ("adjoint functors") going each way between the categories such that each element in a category defines a certain property in the other category and the map carries the element to the self-predicative universal for that property.

$$\Delta: \wp(U) \to \wp(U)\times\wp(U) \text{ and } \cap: \wp(U)\times\wp(U) \to \wp(U)$$
Example of Adjoint Functors Between Partial Orders

The notion of a pair of adjoint functors is ubiquitous; it is one of the main tools that highlights self-predicative universals throughout modern mathematics.

## Heteromorphisms: The Joy of Hets

### Heteromorphisms and Adjunctions

We have seen that there are two self-predicative universals (often one is trivial like $\Delta(c)$ in the above example) involved in an adjunction and that the object-to-object maps or relations were always within one category (or partial order), e.g., the "hom-sets" in a category where "hom" is short for homomorphism (a morphism between objects in the *same* category). Using object-to-object maps between objects of *different* categories (properly called "heteromorphisms" or "chimera morphisms"), the notion of an adjunction can be factored into two semi- or half-adjunctions, each of which isolates a self-predicative universal [Ellerman 2006, 2007].

---

[9] We have written $\langle a,b\rangle \supseteq \Delta(c)$ "backwards" to be an instance of the participation relation "x μ $u_F$" but it would more conventionally be written $\Delta(c) \subseteq \langle a,b\rangle$.

This heteromorphic treatment of adjoints will be illustrated using the above example. The objects $c \in \wp(U)$ in the partial order $\wp(U)$ are single subsets c of U and the objects <a,b> in the partial order $\wp(U) \times \wp(U)$ are pairs of subsets of U. A *heteromorphism* or *het* from a single subset c to the pair of subsets <a,b> is given by the "fork" $c \subseteq a$ and $c \subseteq b$ which could be symbolized c → <a,b>. Fixing <a,b>, there is a single subset $\cap(<a,b>) = a \cap b$ with a canonical het $a \cap b$ → <a,b>. Then the functor that takes <a,b> to $a \cap b$ gives a *right semi-adjunction* if for every het from c → <a,b>, there is a (unique) hom $c \subseteq a \cap b$ that gives us the following (commutative) diagram where the arrows are hets.

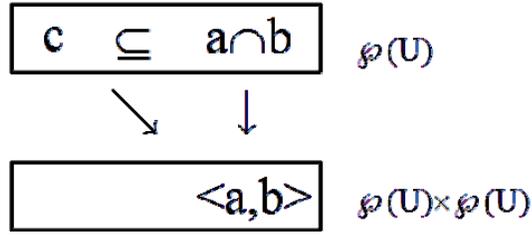

Fig. 1: Right semi-adjunction diagram.

Hence the canonical het $a \cap b$ → <a,b> is not only a universal for the property c → <a,b> but is self-predicative by the identity morphism. This gives the following if-and-only-if equivalence between the diagonal het and the horizontal hom:

$$c \to <a,b> \text{ iff } c \subseteq a \cap b.$$

Universality condition for right semi-adjunction.

Similarly for given c, there is a canonical het c → Δ(c) given by the two identity maps. Then the functor that takes c to Δ(c) is a *left semi-adjunction* if for any given het c → <a,b>, there is a (unique) Δ(c) ⊆ <a,b> to make the following diagram commute.

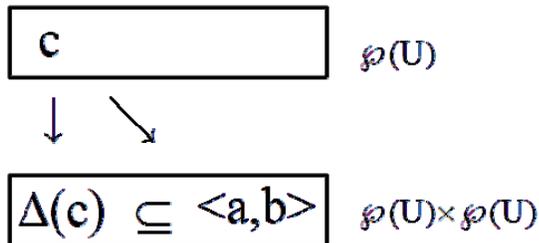

Fig. 2: Left semi-adjunction diagram.

Hence the canonical het c → Δ(c) is not only a universal for the property c → <a,b> but is self-predicative by the identity morphism. The corresponding universality equivalence is:

$$\Delta(c) \subseteq <a,b> \text{ iff } c \to <a,b>$$

Universality condition for left semi-adjunction.

The concept of a semi-adjunction is the most general concept of a self-predicative universal in category theory. The name "semi-adjunction" (or "half-adjunction") is derivative from "adjunction" since two semi-adjunctions with the same diagonal hets, c→<a,b> in this case, combine to give an adjunction:

$$\Delta(c) \subseteq \langle a,b \rangle \text{ iff } c \to \langle a,b \rangle \text{ iff } c \subseteq a \cap b$$

Adjunction equivalence with het middle term.

Gluing together the two left and right semi-adjunction diagrams along the common het c→<a,b> gives the adjunctive square diagram representing an adjunction.

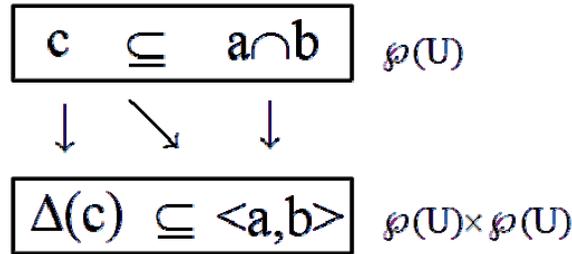

Fig. 3: Two semi-adjunctions = One adjunction

In the orthodox (i.e., "heterophobic" or homs-only) treatment of adjunctions in category theory, the middle het term is left out, so we get the usual form (without any mention of hets) of the adjunction equivalence which is a natural isomorphism in the general case:

$$\Delta(c) \subseteq \langle a,b \rangle \text{ iff } c \subseteq a \cap b$$

Usual form of adjunction equivalence.

The left or right semi-adjunctions are the most general form of self-predicative universals, and the ubiquitous adjunctions are the special cases where left and right semi-adjunctions exist for the same hets.

### Brain functors

If the self-predicative universals of category theory, which combine in one way to form an adjunction, serve to delineate the important, paradigmatic, canonical, or essential concepts and structures within pure mathematics, then one might well expect the self-predicative universals to also be important in applications.

By explicitly adding heteromorphisms to the usual homs-only presentation of category theory, the theory can directly represent interactions between the objects of different categories—which is why hets have long been a part of mathematical practice. For instance, before the homs-only (or "heterophobic") treatment of the free group construction, the universal mapping property of the free group was stated using *both* hets (→) and homs (⇒) :

> for any set X, there is a group F[X] and a set-to-group map ("het") X→F[X] (injection of the generators) such that for any other set-to-group map ("het") X→G, there is the group homomorphism F[X] ⇒G such that
>
> X→G = X→F[X] ⇒G.

Thus the traditional statement of the free group property simply says that the free group functor gives a left semi-adjunction (and does not mention any right semi-adjunction). In the homs-only presentation, the "device" of the underlying set functor U is used as a right semi-adjunction to state (in a rather stilted manner) essentially the same fact of the universal mapping property of free groups without mentioning hets.

for any set X, there is a group F[X] and a set homomorphism X⇒U[F[X]] ("injection of the generators") such that for any other set homomorphism X⇒U[G] to the underlying set of a group G, there is the group homomorphism F[X] ⇒G such that for the image of the underlying set functor, U[F[X]] ⇒U[G],

$$X \Rightarrow U[G] = X \Rightarrow U[F[X]] \Rightarrow U[G].$$

In many adjunctions, the *important* fact is expressed by either the left or right semi-adjunction, but the "device" of the other semi-adjunction is used to express the fact in a het-free manner [Ellerman 2006].

Another payoff from analyzing the important concept of an adjunction into two semi-adjunctions (in addition to the device-free presentation of the important part of an adjunction) is that we can then reassemble those parts in a different way to define the cognate concept speculatively named a "brain functor." The basic intuition is to think of one category in a semi-adjunction as the "environment" and the other category as an "organism." Instead of semi-adjunctions representing within each category the hets going one way between the categories, suppose the hets going both ways were represented by semi-adjunctions within one of the categories (the "organism").

A het from the environment to the organism represents, say, a visual or auditory stimulus. Then a left semi-adjunction would play the role of the brain in providing the re-cognition (expressed by the intentionality-of-perception slogan: "seeing is seeing-as") of the stimulus as a perception of, say, a tree where the internal re-cognition is represented by the morphism ⇒ inside the "organism" category.

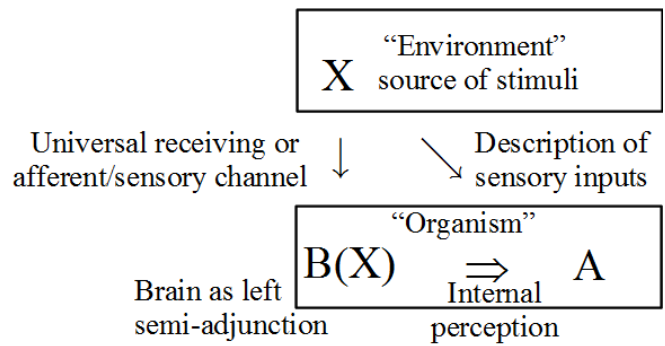

Fig. 4: Perceiving brain presented as a left semi-adjunction.

Perhaps not surprisingly, this mathematically models the old philosophical theme in the Platonic tradition that external stimuli do not give knowledge; the stimuli only trigger the internal perception, recognition, or recollection (as in Plato's *Meno*) that is knowledge. In *De Magistro* (The Teacher), the neo-Platonic Christian philosopher Augustine of Hippo developed an argument (in the form of a dialogue with his son Adeodatus) that as teachers teach, it is only the student's internal appropriation of what is taught that gives understanding.

Then those who are called pupils consider within themselves whether what has been explained has been said truly; looking of course to that interior truth, according to the measure of which each is able. Thus they learn,…. But men are mistaken, so that they call those teachers who are not, merely because for the most part there is no delay between the time of speaking and the time of cognition. And since after the speaker has reminded them, the pupils quickly learn within, they think that they have been taught outwardly by him who prompts them. (Augustine *De Magistro*, Chapter XIV)

The basic point is the active role of the mind in *generating* understanding. This is clear even at the simple level of understanding spoken words. We hear the auditory sense data of words in a completely strange language as well as the words in our native language. But the strange words, like @#$%^, bounce off our minds with no resultant understanding while the words in a familiar language prompt an internal process of generating a meaning so that we *understand* the words. Thus it could be said that "understanding a language" means there is a left semi-adjunction for the heard statements in that language, but there is no such internal re-cognition mechanism for the heard auditory inputs in a strange language.

There are also hets going the other way from the "organism" to the "environment" and there is a similar distinction between mere behavior and an action that expresses an intention. Mathematically that is described by dualizing or turning the arrows around which gives an acting brain presented as a right semi-adjunction.

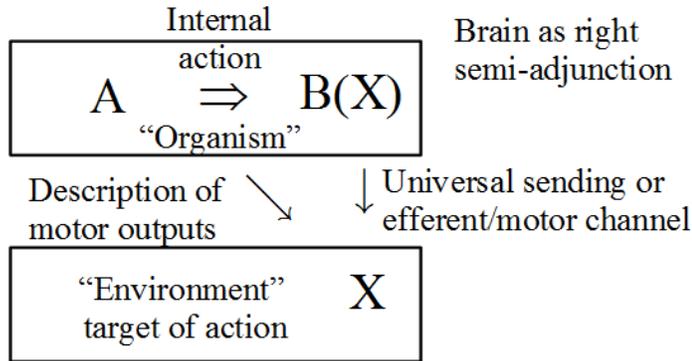

Fig. 5: Acting brain as a right semi-adjunction.

In the heteromorphic treatment of adjunctions, an adjunction arises when the hets from one category to another, Het(X,A), have a right semi-adjunction, Het(X,A) ≅ Hom(X,G(A)), *and* a left semi-adjunction, Hom(F(X),A) ≅ Het(X,A). But instead of taking the same set of hets as being represented by two different functors on the right and left, suppose we consider a single functor B(X) that represents the hets Het(X,A) on the left:

$$\text{Het}(X,A) \cong \text{Hom}(B(X),A),$$

*and* represents the hets Het(A,X) [going in the opposite direction] on the right:

$$\text{Hom}(A,B(X)) \cong \text{Het}(A,X).$$

If the hets each way between two categories are represented by the same functor B(X) as left and right semi-adjunctions, then that functor is said to be a *brain functor*. Thus instead of a pair of functors being adjoint, we have a single functor B(X) with values within one of the categories (the "organism") as representing the two-way interactions, "cognition" and "action," between that category and another one (the "environment").

**Why hets?**

One standard type of brain functor is provided by any functor with both a left and right adjoint (both of which can be described in the usual het-free manner).[10] But the concept of a brain functor is defined using the notion of heteromorphisms as morphisms between the objects of different categories. The use of hets allows a new range of applications of category theory since it allows these "chimera morphisms" to model the direct influence of objects of one category (e.g., "environment") on objects of another category (e.g., "organism")—and vice versa.

An adjunction can be viewed as one way of putting together the building blocks of semi-adjunctions, and the concept of a brain functor is the cognate concept obtained by putting the building blocks of semi-adjunctions together in another way. The adjunctive diagram for an adjunction arises by gluing together the two diagrams for the left and right semi-adjunctions along the common diagonal het X→A.

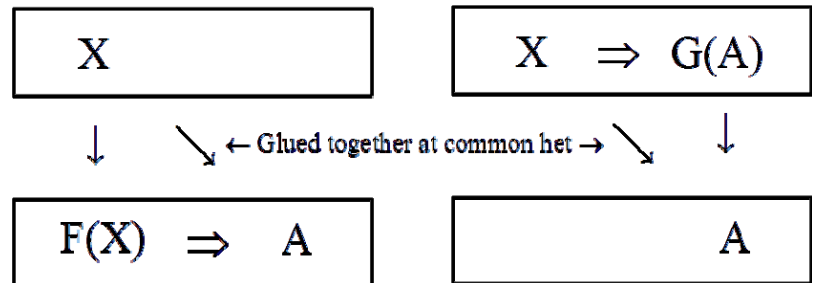

Fig. 6: Combining two semi-adjunctions make an adjunction

The combined diagram is the *adjunctive square* that represents an adjunction in the hetermorphic treatment of adjunctions.

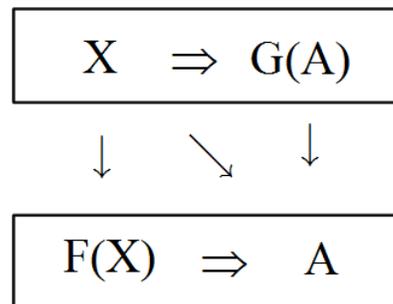

Fig. 7: Adjunctive Square Diagram.

The adjunctive square diagram might be compared to the usual adjunction diagram in the het-free treatment which shows only homs ($\Rightarrow$) inside each category (represented by rectangles) with no hets ($\rightarrow$) between the objects of different categories.

---

[10] The aforementioned underlying set functor that takes a group G to its underlying set U[G] is a rather trivial example of a brain functor that does not arise from having both a left and right adjoint. It does have a left adjoint (the free group functor) but no right adjoint.

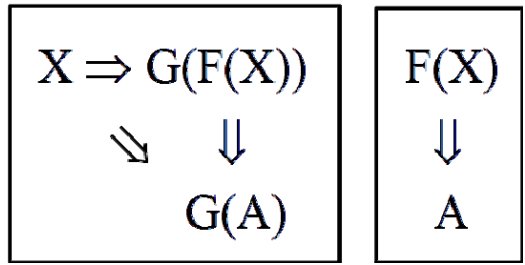
Fig. 8: Usual het-free or homs-only diagram for an adjunction.

The diagram for a brain functor is obtained by gluing together the diagrams for the left and right semi-adjunctions at the common values of the brain functor B(X).

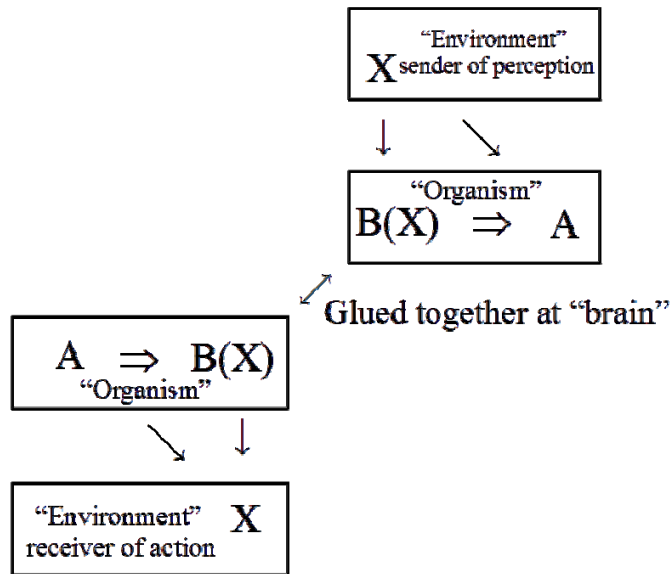
Fig. 9: Combining two semi-adjunctions to make a "brain"

This gives the brain functor "butterfly"[11] diagram—where we have taken the liberty to relabel the diagram for the brain as the language faculty for understanding and producing speech.

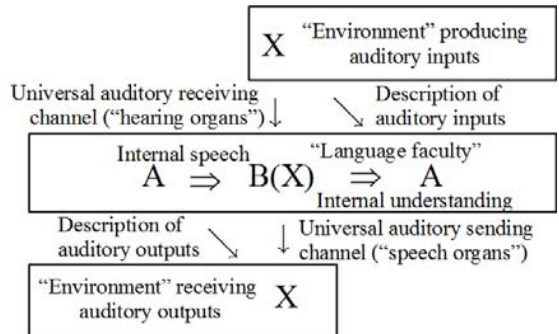
Fig. 10: Brain functor diagram interpreted as language faculty

---

[11] Erase the labels, and then the arrows give the outline of a butterfly.

Wilhelm von Humboldt recognized the symmetry between the speaker and listener, which in the same person is abstractly represented as the dual functions of the "selfsame power" of the language faculty in the above butterfly diagram.

> Nothing can be present in the mind (Seele) that has not originated from one's own activity. Moreover understanding and speaking are but different effects of the selfsame power of speech. Speaking is never comparable to the transmission of mere matter (Stoff). In the person comprehending as well as in the speaker, the subject matter must be developed by the individual's own innate power. What the listener receives is merely the harmonious vocal stimulus. [Humboldt 1997, 102]

If both the triangular "wings" could be filled-out as adjunctive squares, then the brain functor would have left and right adjoints. A simple example of such a brain functor, where the two categories are partial orders, is given in the Appendix—along with a more complex conceptual model.

**Category Theory and Foundations**

Over the last half-century, there has been some antagonism between category theory and set theory on the "foundations question." But we have seen that the two theories complement each other as opposite bookends qua theories of universals (always-self-predicative or never-self-predicative universals). Hence there are some grounds for peaceful coexistence since both theories have distinct foundational roles as theories of universals.

There is also controversy about how category theory is relevant to the foundations of mathematics. The safest view is that category theory is of foundational importance in several different ways that are not mutually exclusive. One approach is to provide one overarching theory in which all or most mathematics can be formulated—as in the original approach of set-theoretic foundations. No one has been more associated with category theory and foundations than William Lawvere, and he seems to have emphasized several approaches. At first, Lawvere [1966] proposed the *category of categories* as an overarching foundational theory but later Lawvere and Tierney's *theory of topoi* [e.g., Lawvere 1972] played a similar role. Independently, the *Univalent Foundations Program* [2013] has been recently proposed as a foundational theory that is heavily based on category theory.

Another approach is that of category-theoretic structuralism [Awodey 1996; Landry 2006]. The morphism language of category theory is the natural language to isolate and describe the structure that objects have independent of the nature of the substratum or other particulars of the case.

Yet another approach (and central to the topic of this paper) started with the characterization by MacLane [1948] and Samuel [1948] of many standard and important structures in mathematics using universal mapping properties.[12] Note that this isolation of the concept of the universals came after Eilenberg and MacLane's [1945] introduction of categories, functors, and natural transformations, and it took another decade before Kan's [1958] isolation of the concept of adjoint functors. Lawvere [1969] has also considered this *universals-approach* where the foundational importance of category theory is that it provides the tools to characterize "what is universal in mathematics" [Lawvere 1969, 281]—"the assumption being that what is universal is to be revealed by adjoint functors." [Landry and Marquis 2005, 12] That is essentially the approach taken here [and in Ellerman 1988, 1995]. Category theory's foundational relevance is that it provides the concepts of universality to characterize the *important* structures, schema, or Forms throughout mathematics.

---

[12] It is interesting to see how MacLane handled the hets naturally involved in a universal mapping property before the development of the homs-only devices used later. To characterize the product within a category C, a het $c \to \langle a,b \rangle$ from an object c in the category C to an object $\langle a,b \rangle$ in the category C×C is a pair of C-morphisms $f: c \Rightarrow a$ and $g: c \Rightarrow b$ which MacLane called a "system" of maps [1948, 264] and was later called a "cone" (or "fork"). Fixing $\langle a,b \rangle$, the product a×b is given by what is here called a right semi-adjunction. That is, there is a universal het a×b $\to \langle a,b \rangle$ consisting of the cone-het of C-homs $p_a: a \times b \Rightarrow a$ and $p_b: a \times b \Rightarrow b$ (the projections) from a×b to $\langle a,b \rangle$ such that for any cone-het $f: c \Rightarrow a$ and $g: c \Rightarrow b$ from an object c in C to the object $\langle a,b \rangle$ in C×C, there is a unique C-hom

Moreover, if the self-predicative universals of category theory characterize important concepts within pure mathematics, it should not be too surprising if they might also characterize, albeit at a very abstract level, important concepts and canonical schema (or Forms) in applications. The application scheme outlined here is the brain functor (obtained by rearranging the heteromorphic building blocks of adjunctions) which abstractly models the dual functions of perception and action.

The importance of category theory in mathematics is that it provides a criterion of importance in mathematics. Category theory provides the concepts to isolate the universal instance (where it exists) from among all the instances of a property. The Self-Predicative or Concrete Universal is the most important instance of a property because it represents the property in a paradigmatic way. All instances have the property by virtue of participating in the Self-Predicative Universal.

**Appendix: Two examples of brain functors**

A simple example of a brain functor using partial orders will be developed first. In this setting, only the simplest "brain function" can be modeled, namely the building and functioning of an internal model of the external reality such as an internal coordinate system to map an external set of locations. The external reality is given by a set of atomic points or locations Y, the atomic coordinates are the points in X, and the coordinate mapping function is the given function f:X→Y. Just to keep the mathematics not completely trivial, we do not require f to be an isomorphism; multiple coordinates might refer to the same point (i.e., f is not necessarily one-to-one) and some points might not have coordinates (i.e., f is not necessarily onto). The two partial orders are the inclusion-ordered subsets of points $\wp(Y)$ and the inclusion-ordered subsets of coordinates $\wp(X)$.

---

c⇒a×b such that the C-hom followed by the universal het equals the given het. Thus the C-hom c⇒a×b is the unique ("participation") factor map that "transfers" (in this case "reflects") the property of "*being a het from any c in C to the fixed <a,b> in C×C*" from the universal such cone-het (the projections) to the given such cone-het f:c⇒a and g:c⇒b. That is the het statement of the universal mapping property for the product without the "device" of the diagonal functor C→C×C used in the homs-only treatment [e.g., MacLane 1971, 68]. Samuel worked with structured sets such as S-sets and T-sets so the homs would be the S-mappings or T-mappings and were denoted by ordinary alphabetical letters. He called the hets from an S-set to a T-set "(S-T)-mappings" [Samuel 1948, 592] and he denoted them with Greek letters to recognize their difference from homs. Samuel even noted that the composition of a (S-T)-mapping with a T-mapping would be a (S-T)-mapping, i.e., hets composed with homs are hets [see Ellerman 2006 for a general treatment].

In the case where the "brain" is an "electronic brain" or computer, Y is the set of locations on an external input/output device such as a floppy disk or any other external memory device. Each location is marked with a 0 or 1, so the subsets $V \in \wp(Y)$ would be the external sets of 1s. The set X would be the set of internal memory locations which also contain either a 0 or 1, so the subsets $U \in \wp(X)$ are the internal sets of 1s. The coordinate function $f: X \to Y$ maps the internal memory locations to the external disk locations. The dual perception/action functions in the electronic brain would be the familiar read/write operations between the computer and the external input/output device.[13]

The brain functor in this example is $f^{-1}: \wp(Y) \to \wp(X)$ where for any subset $V \in \wp(Y)$, the value of the brain functor is:

$$f^{-1}(V) = \{x \in X : f(x) \in V\}.$$

Given a subset $V \subseteq Y$, what is the "best" internal subset $U \subseteq X$ that represents or recognizes V? The heteromorphism $V \to U$ is defined by the property, $F(U) =$ "U is *complete for V*" in the sense that all the $x \in X$ that map to V are contained in U, i.e.,

$$V \to U \text{ means } V \subseteq \{y \in Y : \forall x, \text{ if } f(x) = y \text{ then } x \in U\}.$$

The left semi-adjunction for the property $F(U) =$ "U is complete for V" is given by the smallest complete subset $f^{-1}(V) \in \wp(X)$ and the universality condition: $f^{-1}(V) \subseteq U$ iff $V \to U$, is satisfied.

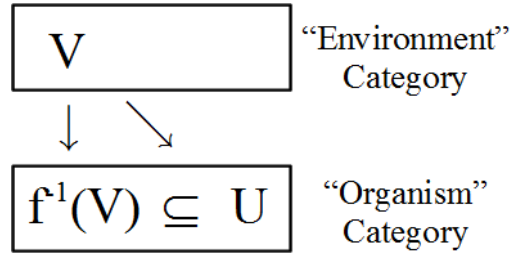

Fig. 11: Left semi-adjunction to "read" V
with smallest complete subset $f^{-1}(V)$

The self-predicative universal $f^{-1}(V) \in \wp(X)$ has the property, i.e., $V \to f^{-1}(V)$, and a subset $U \in \wp(X)$ has the property, i.e., $V \to U$, if and only U participates in the self-predicative universal $U \supseteq f^{-1}(V) = u_F$ (where "participation" is written as the reverse inclusion).

In the dual case of "action," the het $U' \to V$ going in the opposite direction from a subset U' of X to a subset V of Y is defined by the property, $G(U') =$ "U' is *consistent with V*" in the sense that no coordinate in U' maps outside of V, i.e.,

$$U' \to V \text{ means } f(U') \subseteq V.$$

The right semi-adjunction for the property $G(U') =$ "U' is consistent with V" is given by the largest consistent subset $f^{-1}(V) \in \wp(X)$, and the universality condition is: $U' \subseteq f^{-1}(V)$ iff $U' \to V$.

---

[13] See also the coordinate-plot scheme in Lawvere and Schanuel 1997, 86.

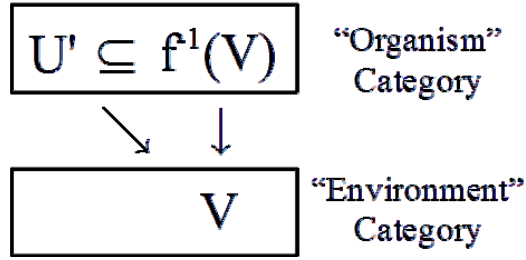

Fig. 12: Right semi-adjunction to "write" V
with largest consistent subset $f^{-1}(V)$.

The self-predicative universal $f^{-1}(V) \in \wp(X)$ has the property, $f^{-1}(V) \to V$, and a subset $U' \in \wp(X)$ has the property, i.e., $U' \to V$, if and only if $U'$ participates in the self-predicative universal, i.e., $U' \subseteq f^{-1}(V) = u_G$ (where "participation" is the inclusion). Combining the left and right semi-adjunctions at the common "brain" gives the butterfly diagram.

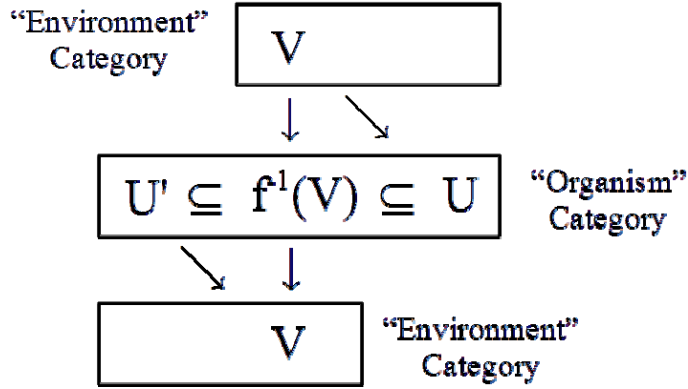

Fig. 13: Butterfly diagram for the brain functor $f^{-1}: \wp(Y) \to \wp(X)$

    Mathematically, this example is an instance of the general result that any functor that has both right and left adjoints is a brain functor. The right adjoint of $f^{-1}(V)$ is usually symbolized as:

$$\forall_f(U) = \{y \in Y : \forall x, \text{ if } f(x) = y \text{ then } x \in U\}$$

with the adjunction equivalence:

$$f^{-1}(V) \subseteq U \text{ iff } V \subseteq \forall_f(U)$$

while the left adjoint of $f^{-1}(V)$ is usually symbolized as:

$$\exists_f(U') = f(U') = \{y \in Y : \exists x \in U', f(x) = y\}$$

with the adjunction equivalence:

$$\exists_f(U') \subseteq V \text{ iff } U' \subseteq f^{-1}(V).$$

The quantifier notation is motivated by the special case where $f$ is the projection $f = p_X : X \times Y \to Y$ so for any binary relation $U \subseteq X \times Y$, then $\exists_f(U) = \{y \in Y : \exists x U(x,y)\}$ and $\forall_f(U) = \{y \in Y : \forall x U(x,y)\}$.

    The fact that the read/write functions of an electronic brain can be modeled by the brain functor $f^{-1}: \wp(Y) \to \wp(X)$ is a simple example of modeling a brain *at a conceptual level*—which says little or nothing about the underlying (electronic) mechanisms. As complexity increases exponentially in animal and human brains, brain functors should be similarly seen as only modeling the brain functions at a conceptual level (e.g., Figure 10 giving the scheme for the

language faculty) while saying nothing about the underlying biological and chemical mechanisms.

A more mathematically complex (beyond our restriction to partial orders) and adequate model of a brain (still at a conceptual level) is provided by the functor taking a finite set of vector spaces $\{V_i\}_{i=1,\ldots,n}$ over the same field (or R-modules over a ring R) to the product $\prod_i V_i$ of the vector spaces. Such a product is also the coproduct $\sum_i V_i$ [Hungerford 1974, 173] and that space may be written as:

$$V_1 \oplus \ldots \oplus V_n \cong \prod_i V_i \cong \sum_i V_i.$$

The het from a set of spaces $\{V_i\}$ to a single space $V$ is a "cocone" of vector space maps $\{V_i \Rightarrow V\}$ and the canonical such het is the set of canonical injections $\{V_i \Rightarrow V_1 \oplus \ldots \oplus V_n\}$ with the "brain" at the point of the cocone. The perception left semi-adjunction then might be taken as conceptually representing the function of the brain as integrating multiple sensory inputs into an interpreted perception.

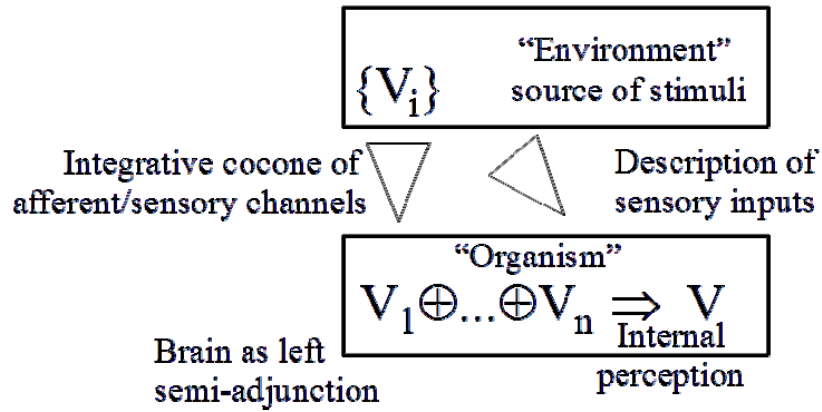

Fig. 14: Brain as integrating sensory inputs into a perception

Dually, a het from single space $V$ to a set of vector spaces $\{V_i\}$ is a cone $\{V \Rightarrow V_i\}$ with the single space at the point of the cone, and the canonical het is the set of projections with the "brain" as the point of the cone: $\{V_1 \oplus \ldots \oplus V_n \Rightarrow V_i\}$. The action right semi-adjunction then might be taken as conceptually representing the function of the brain as integrating or coordinating multiple motor outputs in the performance of an action.

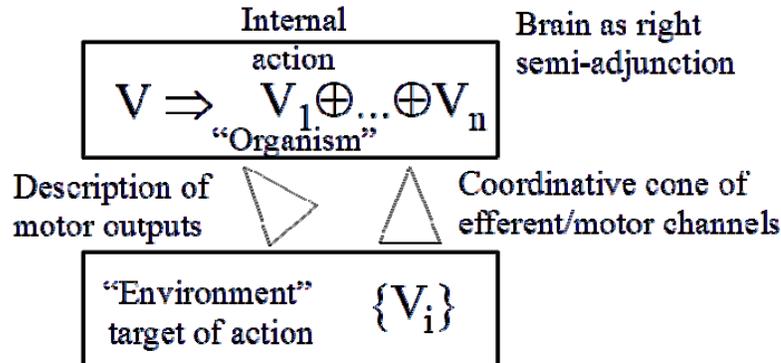

Fig. 15: Brain as coordinating motor outputs into an action

Putting the two semi-adjunctions together gives the butterfly diagram for a brain.

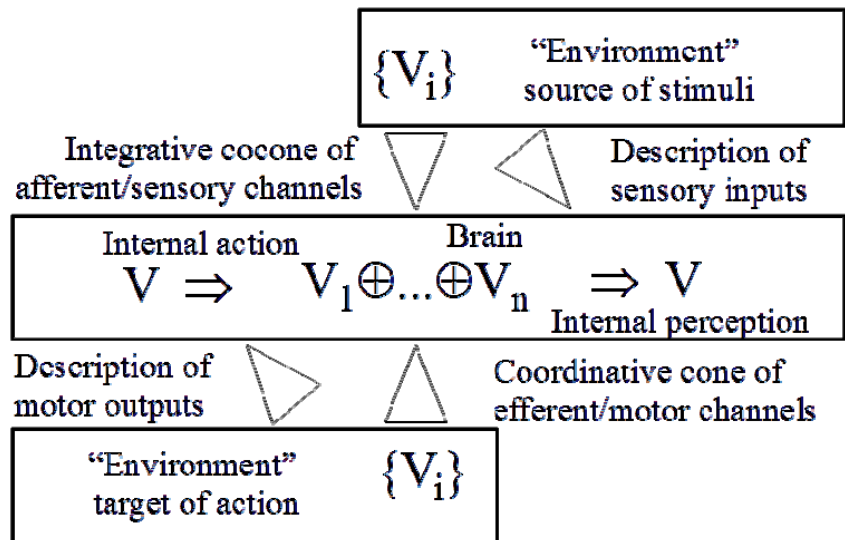

Fig. 15: Conceptual model of a perceiving and acting brain

This gives a conceptual model of a single organ that integrates sensory inputs into a perception and coordinates motor outputs into an action, i.e., a brain.